\documentclass[12pt]{article}
\usepackage{graphics,graphicx,amscd,amssymb,amsmath,amsthm,color}
\usepackage{appendix}
\usepackage{hyperref}

\def\flrna{\lfloor {n\alpha} \rfloor}
\def\flrnb{\lfloor {n\beta} \rfloor}

\def\flria{\lfloor {i\alpha} \rfloor}
\def\flrib{\lfloor {i\beta} \rfloor}
\def\flrnia{\lfloor {(n-i)\alpha} \rfloor}
\def\flrnib{\lfloor {(n-i)\beta} \rfloor}
\theoremstyle{dotless}
\newtheorem{theorem}{Theorem}

\newtheorem{lemma}[theorem]{Lemma}

\theoremstyle{definition}

\newtheorem{remark}[theorem]{Remark}
\newtheorem{definition}[theorem]{Definition}
\newtheorem{notation}[theorem]{Notation}

\newtheorem{example}[theorem]{Example}

\begin{document}
\title{Impartial games whose rulesets produce given continued fractions}
\author{Urban Larsson,\\ Mathematical Sciences \\
Chalmers University of Technology \\ and University of Gothenburg, G\"oteborg, Sweden \\
urban.larsson@chalmers.se\\
\\
Mike Weimerskirch, \\ University of Minnesota \\ weim0024@umn.edu}

\maketitle

\begin{abstract}
We study 2-player impartial games of the form \emph{take-away} which produce P-positions (second player winning positions) corresponding to \emph{complementary Beatty sequences}, given by the continued fractions $(1;k,1,k,1,\ldots)$ and $(k+1;k,1,k,1,\ldots)$. Our problem is the opposite of the main field of research in this area, which is to, given a game, understand its set of P-positions. We are rather given a set of (candidate) P-positions and look for ``simple'' rules. Our rules satisfy two criteria, they are given by a \emph{closed formula} and they are \emph{invariant}, that is, the available \emph{moves} do not depend on the position played from (for all options with non-negative coordinates).
\end{abstract}

\newpage

\section{Introduction}\label{S1}
This paper uses ideas from combinatorial game theory, Beatty sequences, and Sturmian words. We have in many case given the pertinent information in this paper, but have chosen to omit some material on these subjects. The reader who wishes to have more background information on certain topics is directed to the following references:

For standard terminology of impartial removal games on heaps of tokens, see \cite{WW}, for Beatty sequences, see \cite{B}, for $k$-Wythoff Nim, see \cite{W,F}, for Sturmian words, see \cite{L}, for continued fractions, see [K].

Our problem is an inverse to that of the main field of research, which for a given an impartial ruleset $\Gamma$, (for example, $\Gamma = k$-Wythoff Nim) is to determine the P-positions of $\Gamma$ (within reasonable time-complexity). Here we rather start with a particular (candidate) set of P-positions and search for ``simple'' game rules. Let us explain the setting.

Throughout this paper, we will denote the \emph{position} consisting of two heaps of $x\ge 0$ and $y\ge 0$ tokens as $(x,y)$.  When the values of $x$ and $y$ are known, we adopt the convention that $x \le y$, though in general we regard such a position as an unordered multiset, so we identify $(y,x)$ with $(x,y)$. 

Similarly, we let the \emph{move} $(u,v)$ denote a removal of $v>0$ tokens from one of the heaps and $0\le u\le v$ from the other, thus from the position $(x,y)$, the move $(u,v)$ is ambiguous, being either $(x,y)\rightarrow (x-u,y-v)$, provided both $x-u\ge 0$ and $y-v\ge 0$, or $(x,y)\rightarrow (x-v,y-u)$, provided both $x-v\ge 0$ and $y-u\ge 0$. Thus in general, it is necessary to examine both cases. 

Recall that a (homogeneous) \emph{Beatty sequence} is a sequence of integers of the form $(\lfloor n\gamma \rfloor)$, the modulus $\gamma$ being a positive irrational and $n$ ranging over the non-negative integers, here denoted by $\mathbb N$. We are interested in positions of the form 
\begin{align}\label{pairs}
(\flrna,\flrnb)
\end{align} 
for $n \in \mathbb N$, where $0<\alpha < \beta$ are irrationals with
\begin{align}\label{compl}
\alpha ^{-1} + \beta ^{-1} = 1,
\end{align}
that is $1 < \alpha < 2 < \beta $. By (\ref{compl}), the sequences $(\flrna)$ and $(\flrnb)$, where $n$ ranges over the positive integers, $\mathbb Z^+,$ are \emph{complementary}, (see \cite{B}), that is, each positive integer is attained precisely once in precisely one of these sequences.

\begin{figure}[ht!]
  \centering
    \includegraphics[width=0.45\textwidth]{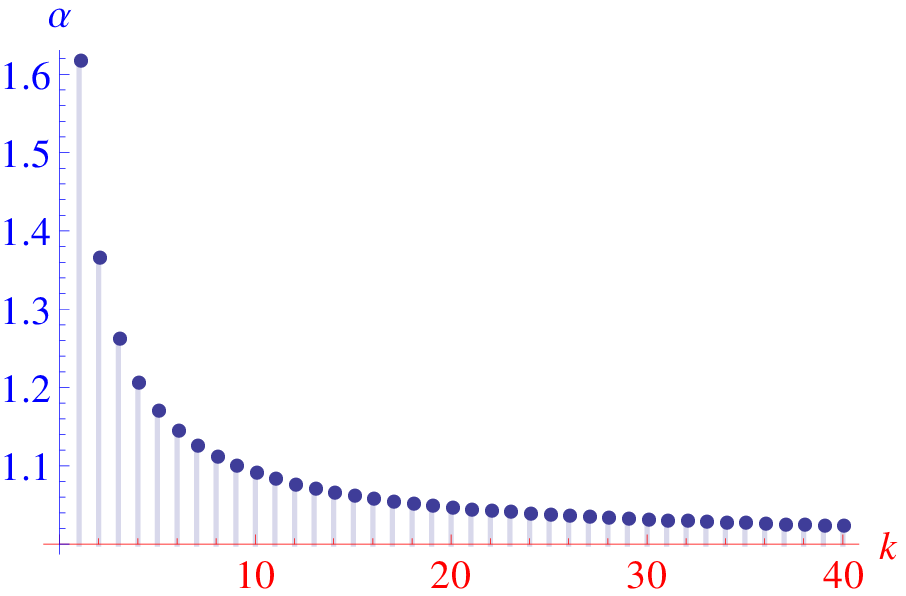} \includegraphics[width=0.45\textwidth]{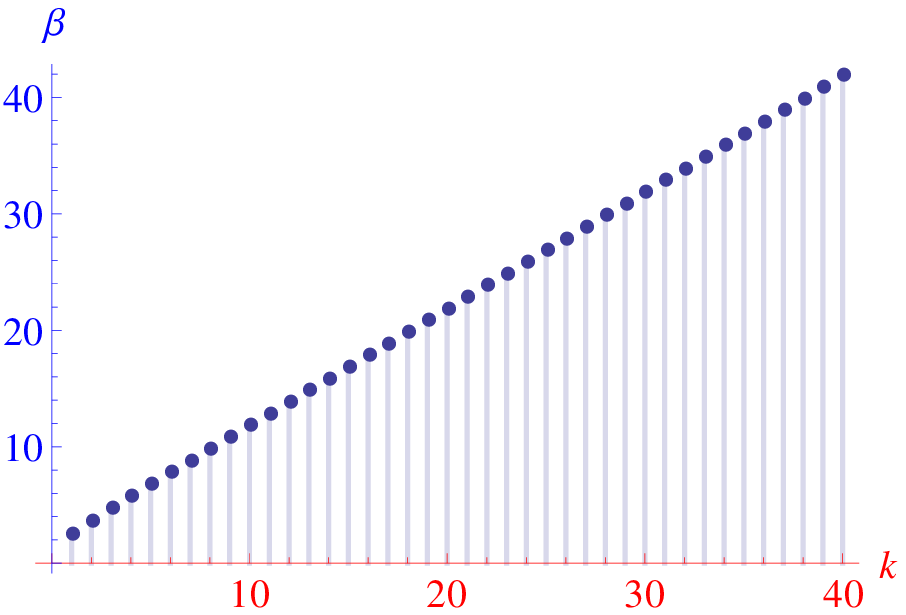} \caption{The numbers $\alpha_k$ and $\beta_k$ for $k\in \{1,40\}$. See also Figure \ref{Ffrac}.}\label{F1}
\end{figure}

In $k$-Wythoff Nim \cite{F}, the P-positions correspond to all unordered pairs of the form in (\ref{pairs}) with 
$$\alpha = [1;k,k,k,\ldots] = \dfrac{2-k+\sqrt{k^2+4}}{2}$$
and 
$$\beta = [k+1;k,k,k,\ldots]=  \dfrac{2+k+\sqrt{k^2+4}}{2}=\alpha + k,$$
where $x = [x_1;x_2,x_3,x_4, \ldots ]$ denotes the unique continued fraction expansion, CF, of $x$, 
\[ 
  x = x_1+\cfrac{1}{x_2+
                                \cfrac{1}{x_3 +
                                \cfrac{1}{x_4 + \cdots.
                                 }}}
\]
For a variation, in \cite{DR} game rules are examined for P-positions corresponding to the CFs 
\begin{align}\label{DR}
[1;1,k,1,k,\ldots] \text{ and } [k+1;k,1,k,1,\ldots]. 
\end{align}
In this paper, we rather study the CF $$\alpha = \alpha_k = [1;k,1,k,1,k,\ldots] = \dfrac{1+\sqrt{1+\frac{4}{k}}}{2}$$ with corresponding $$\beta = \beta_k = [k+1;1,k,1,k,1,\ldots] = \dfrac{k+2+k\sqrt{1+\frac{4}{k}}}{2}=k\alpha + 1=k\alpha^2.$$
Note that $\alpha_k\in (1,1+1/k)$ and $\beta_k\in (k+1,k+2)$. (see Figure \ref{F1} and Lemma \ref{lemma:10})

\begin{notation} For each $k\in \mathbb Z^+$, for all $n\in \mathbb N$, we let $a_n = \flrna$, $b_n = \flrnb$, $c_n=a_n-a_{n-1}$ and $d_n=b_n-b_{n-1}$. Moreover we define the following sequences, $A=(a_1,a_2,\ldots), B=(b_1,b_2,\ldots), C=(c_1,c_2,\ldots)$ and $D=(d_1,d_2,\ldots)$. 
\end{notation}

Then, for all $n\in \mathbb Z^+$, $$b_n = \sum_{j=1}^n d_j\ \text{ and } \ a_n = \sum_{j=1}^n c_j.$$

\begin{example} For $k=2$, 
\begin{align*}
A &= (1,2,4,5,6,8,9,10,12,13,15,16,17,19,20,21,23,24,\ldots )\\
C &= (1,1,2,1,1,2,1,1,2,1,2,1,1,2,1,1,2,1,\ldots)\\
B &= (3,7,11,14,18,22,26,29,33,37,41,44,48,52,55,\ldots )\\
D &= (3,4,4,3,4,4,4,3,4,4,4,3,4,4,3,\ldots)
\end{align*}
\end{example}

Since $\alpha_k \in (1,1+1/k)$ and $\beta_k \in (k+1,k+2)$, for each
$k$ and all $n$, we also get that $$c_n\in \{1,2\} \text{ and } d_n \in \{k+1,k+2\}$$
(see \cite{L}). Moreover, each value is attained infinitely often,
a statement which we strengthen in Section \ref{S2}.

Henceforth, for a fixed $k\in \mathbb Z^+$, let $$S_k=\{(a_n,b_n)\mid n\in \mathbb N \}.$$
Note that the special case of $S_1$ corresponds precisely to the P-positions of Wythoff Nim.

\begin{figure}[ht!]
  \centering
  \includegraphics[width=0.22\textwidth]{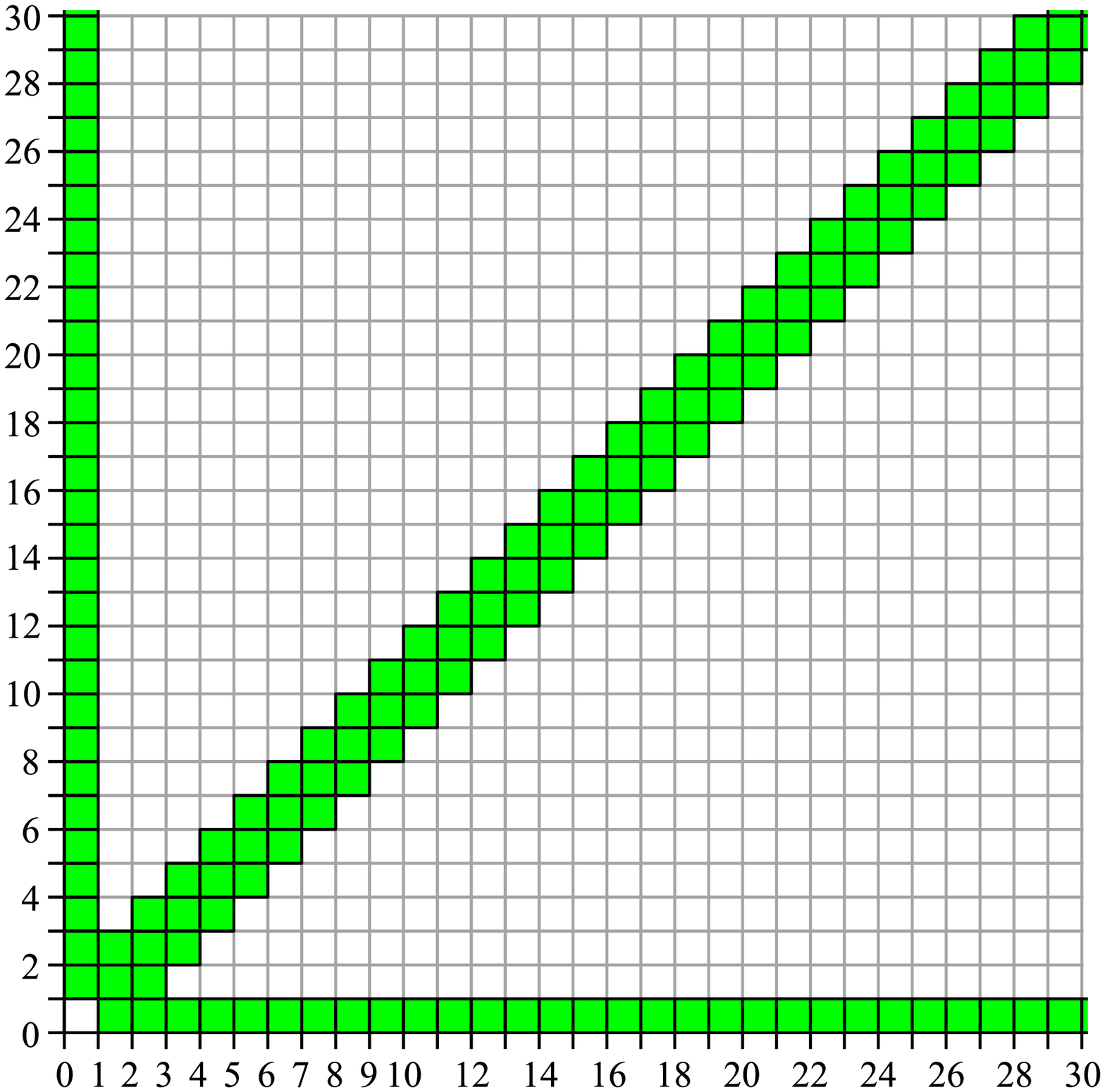} \includegraphics[width=0.22\textwidth]{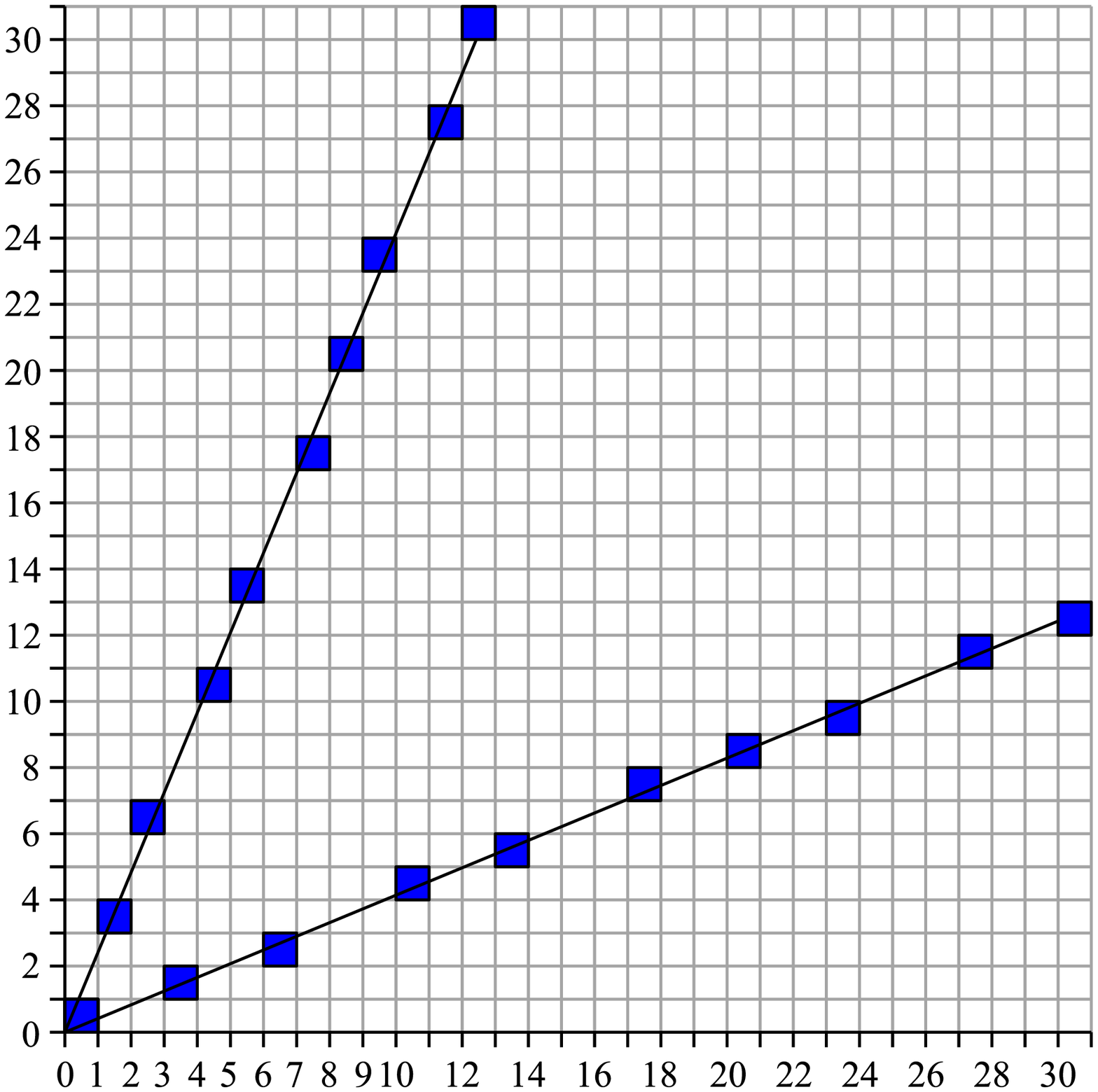}\includegraphics[width=0.22\textwidth]{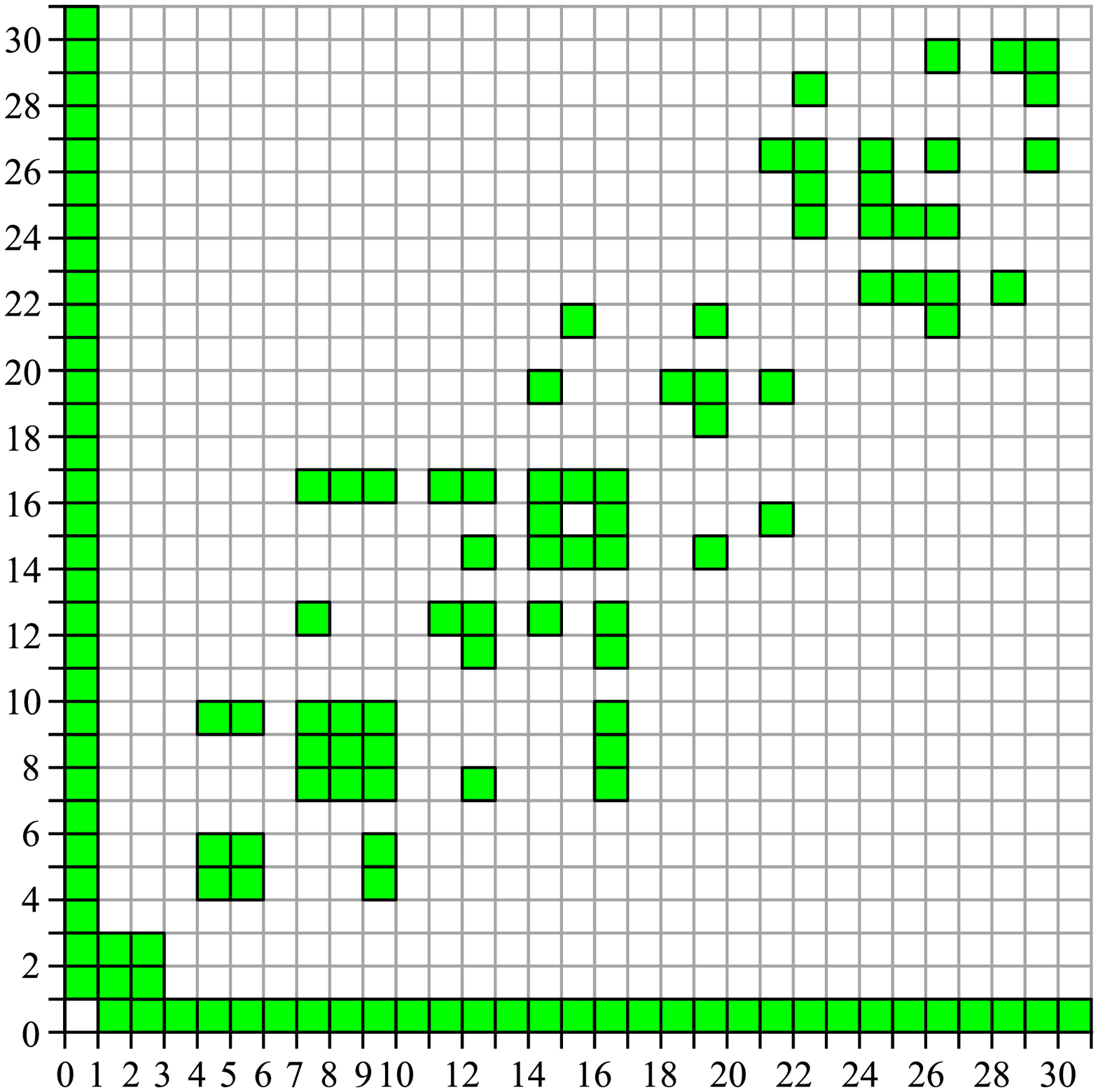} \caption{To the left we display the initial moves of the classical game of $2$-Wythoff Nim and in the middle its initial P-positions together with the corresponding slopes. To the right we give the initial invariant moves for the game $(2\text{-Wythoff Nim})^{\star\star}$, with notation as in \cite{LHF}, which has P-positions of the same form as those of $2$-Wythoff Nim. (The moves are defined via a simple greedy algorithm.)}\label{F2}
\vspace{0.45 cm}
\includegraphics[width=0.22\textwidth]{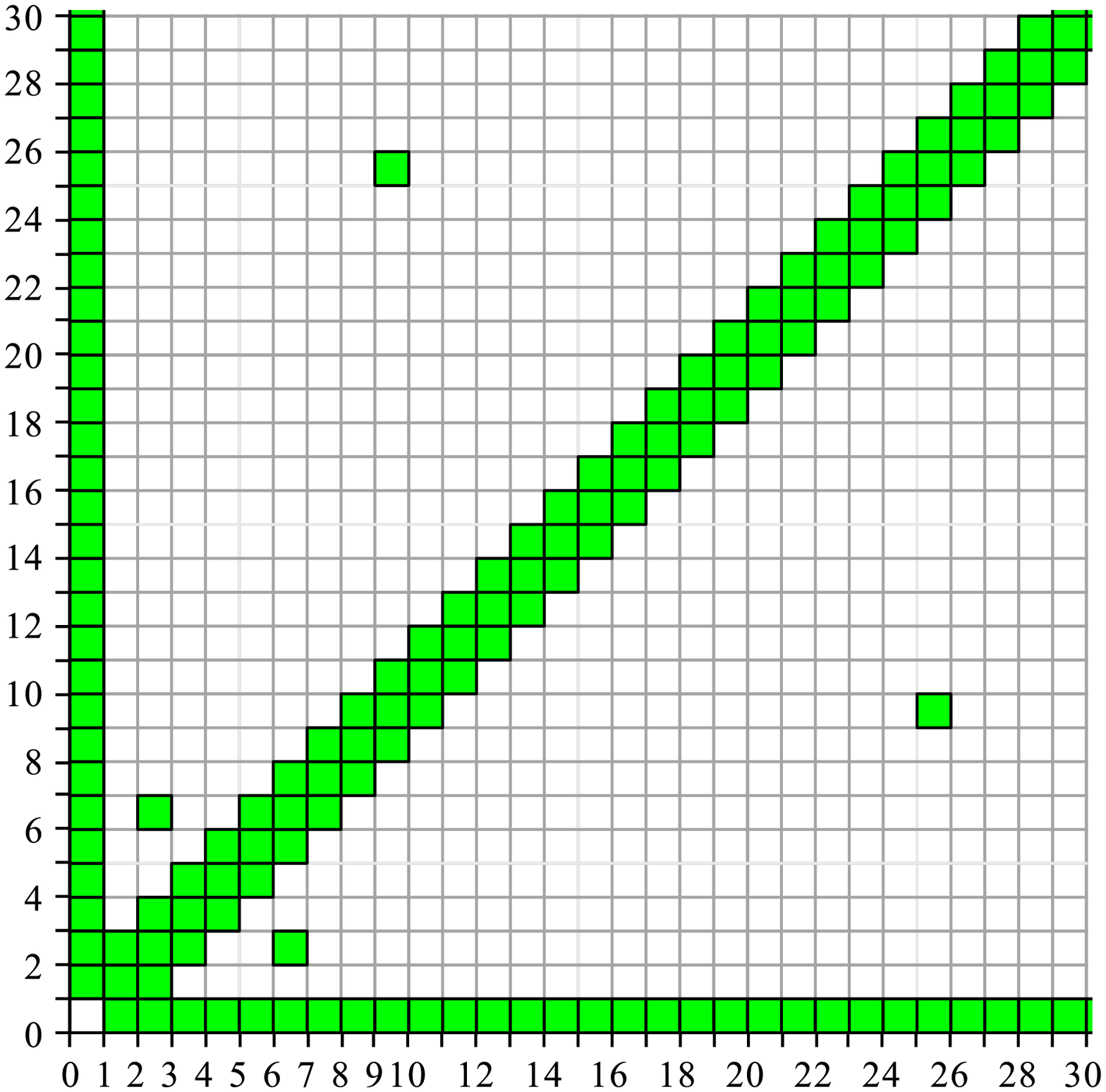} \includegraphics[width=0.22\textwidth]{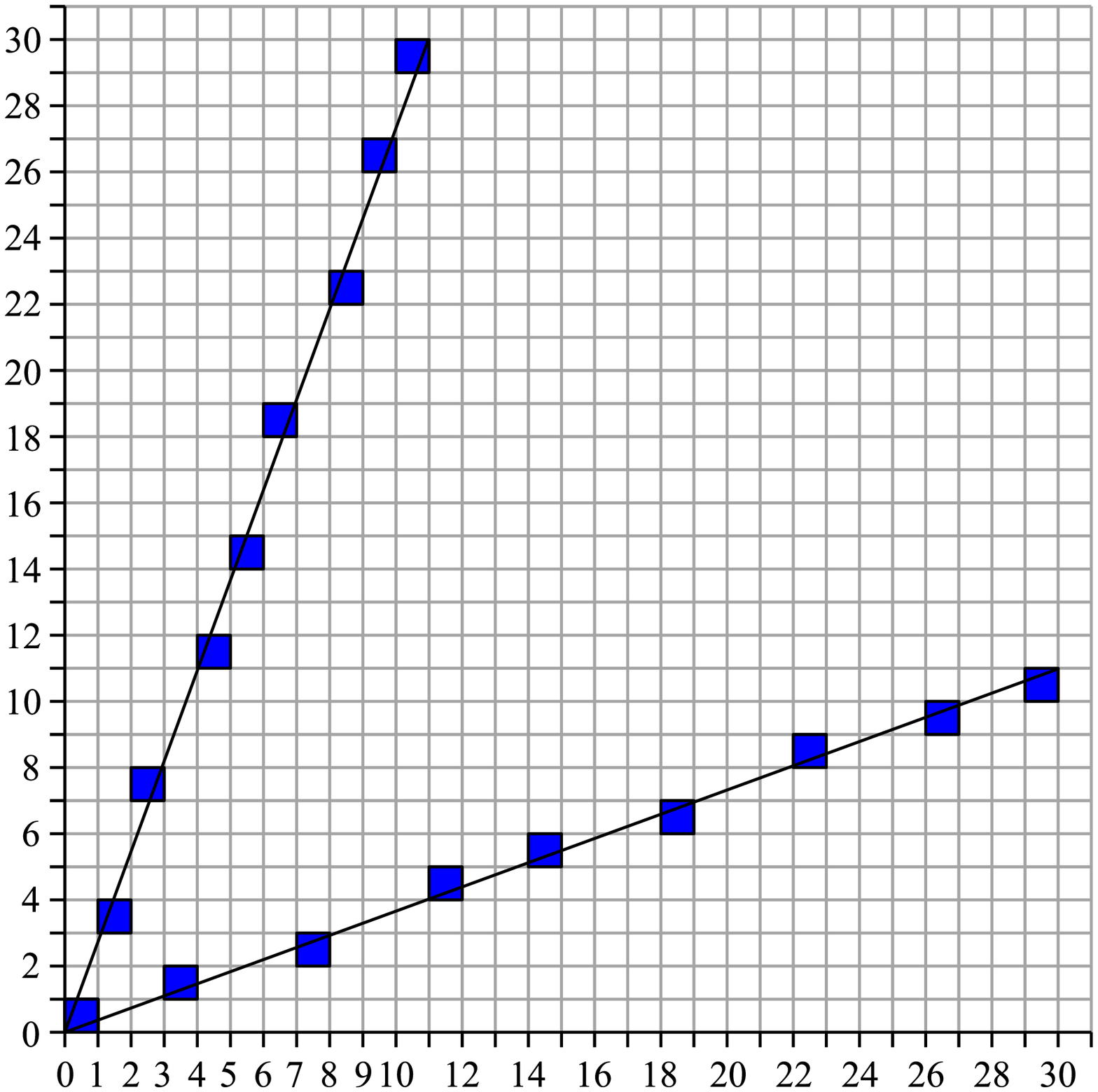}\includegraphics[width=0.22\textwidth]{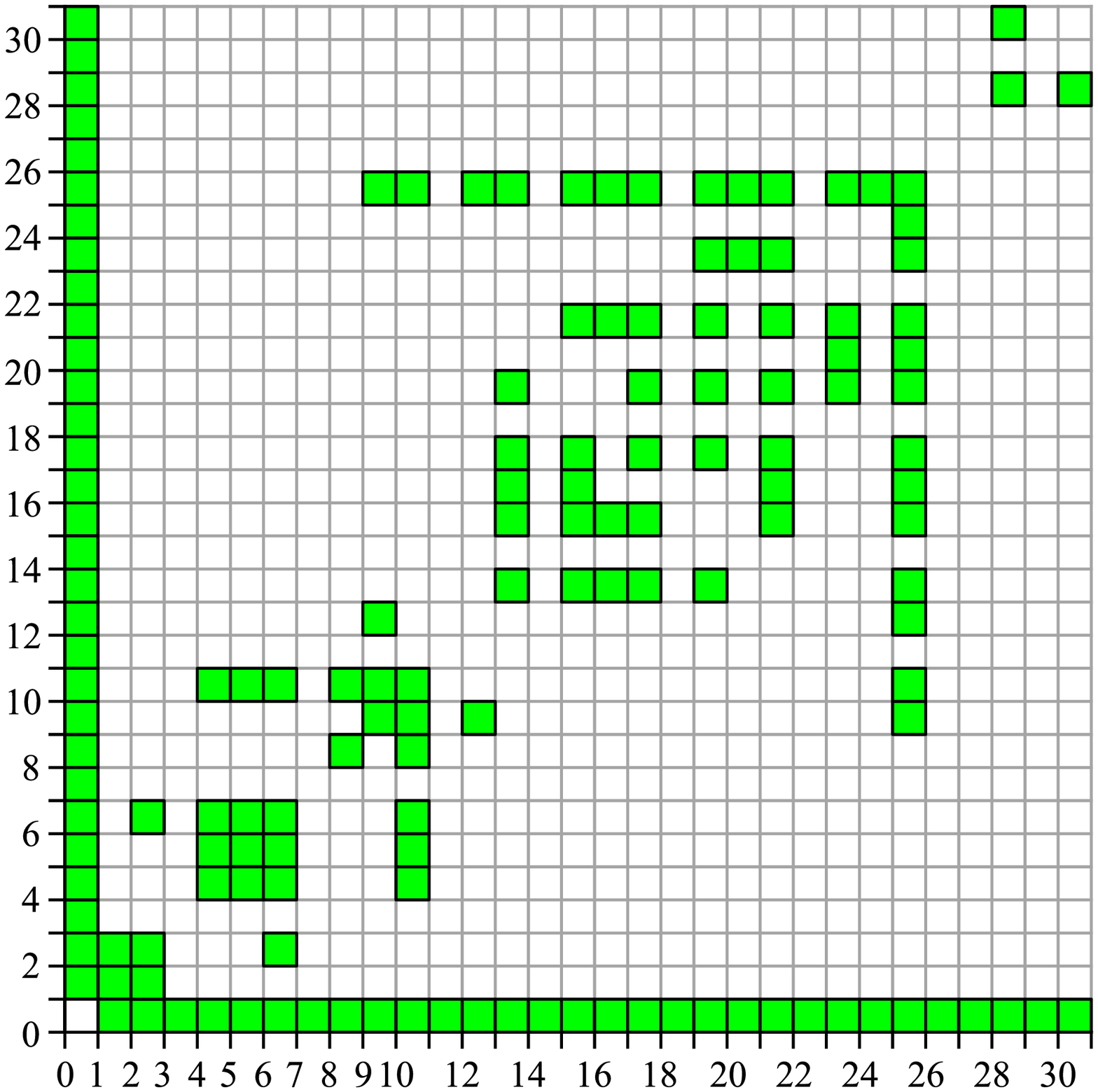}\caption{The left most figure displays the initial moves of our game for $k=2$ as given in Theorem \ref{T1} (see also Example \ref{E4} for the extra moves). In the middle we see the P-positions and to the right the initial moves of the invariant game from \cite{LHF} with P-positions identical to the set $S_2$.}\label{F3}
\vspace{0.45 cm}
  \includegraphics[width=0.22\textwidth]{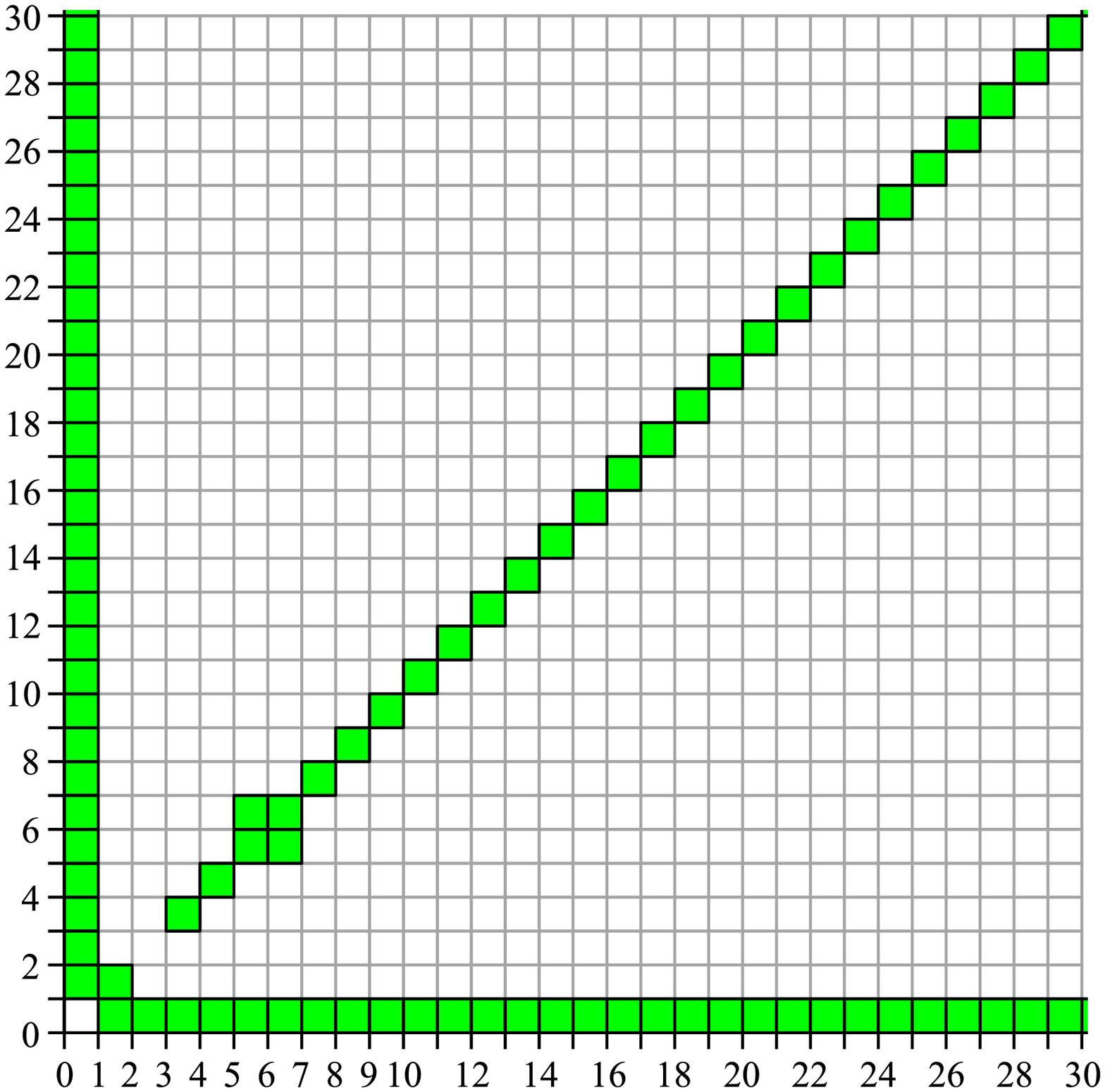} \includegraphics[width=0.22\textwidth]{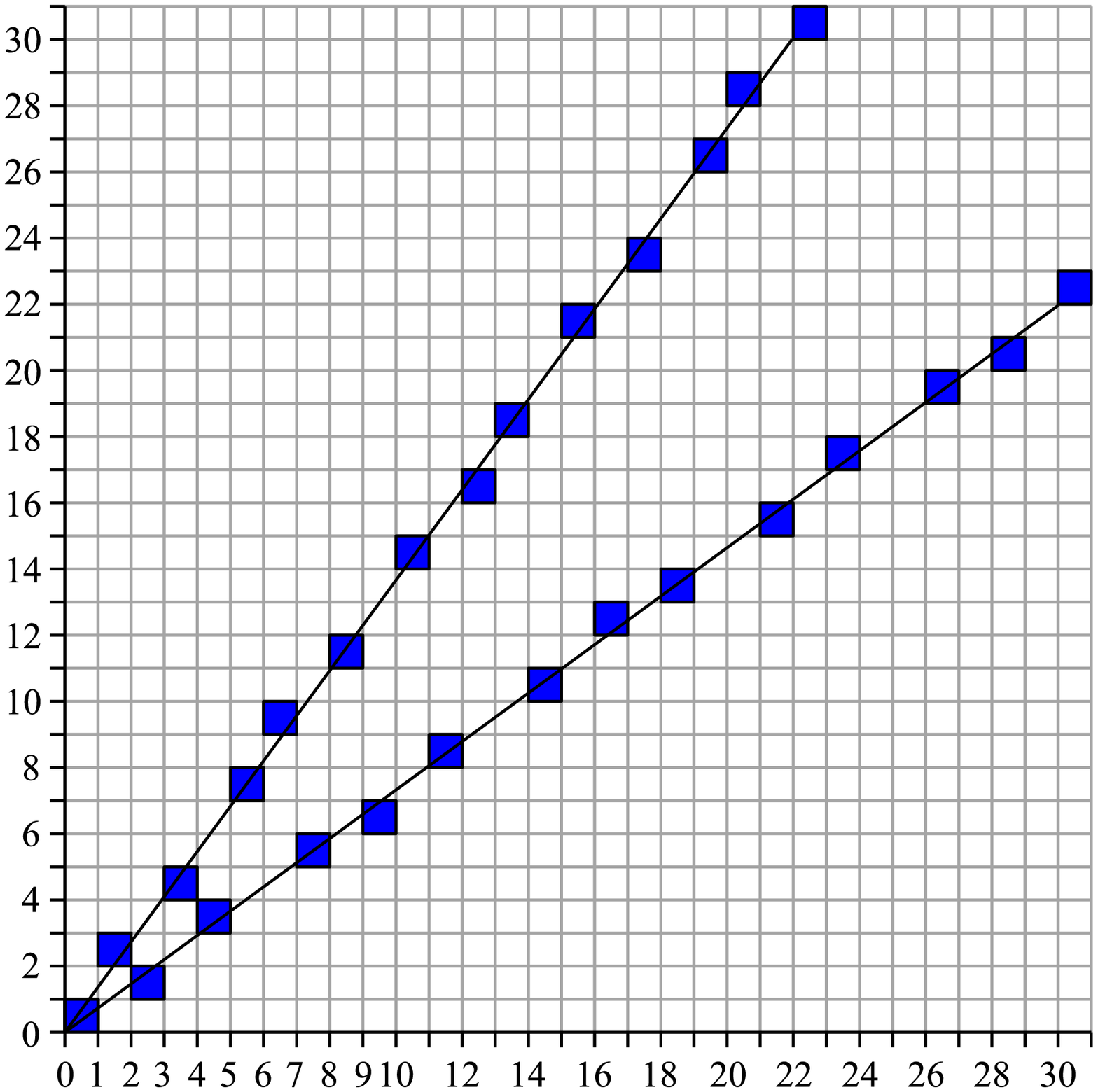}\includegraphics[width=0.22\textwidth]{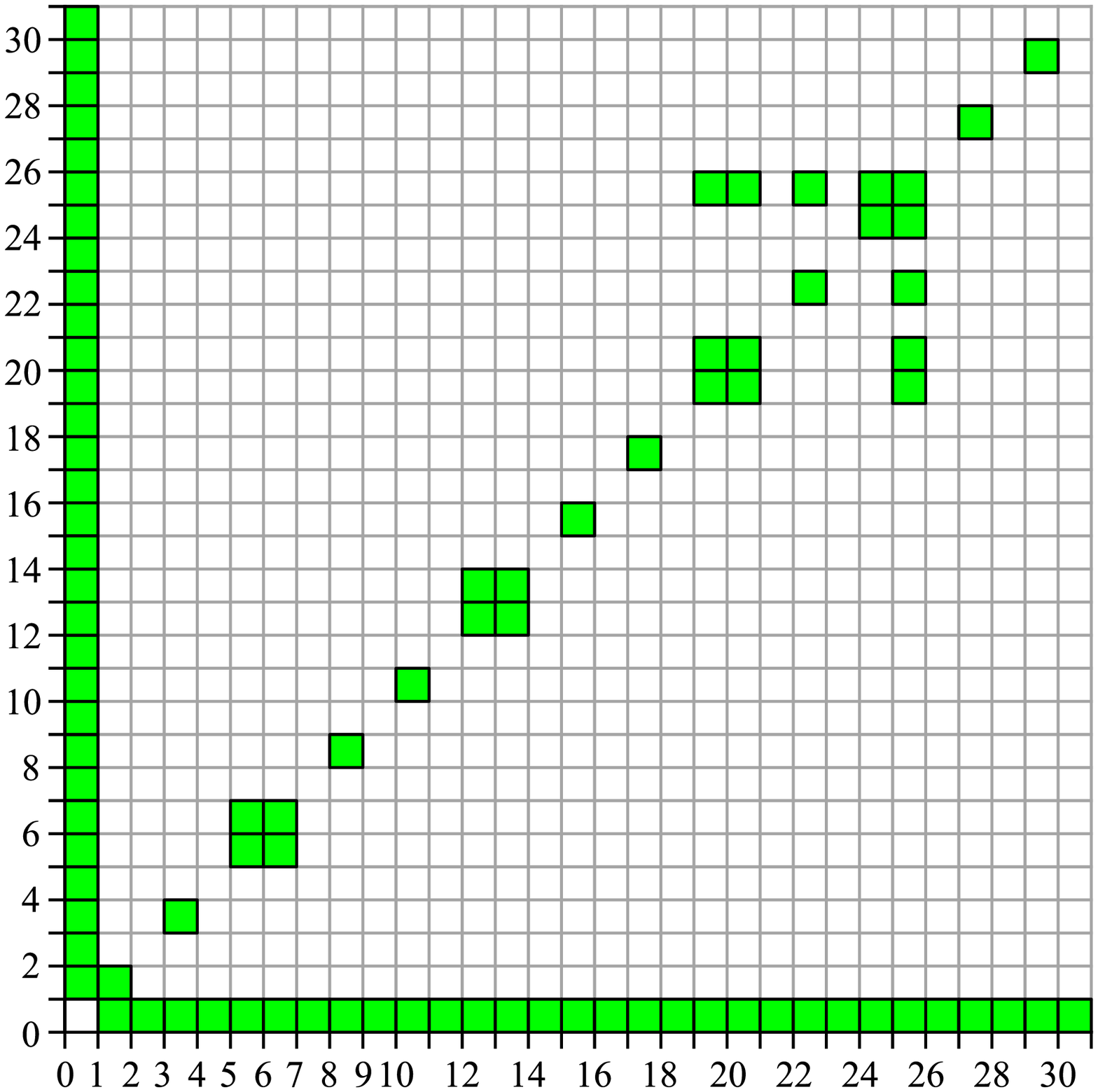} \caption{These figures represent the moves (left figure) resolving the P-positions (middle figure) given by the continued fraction (\ref{DR}) with $k=2$ from \cite{DR}. The right most figure gives the initial moves of the invariant game from \cite{LHF} with identical P-positions.}\label{F4}
\end{figure}

The problem of finding a closed formula ruleset such that the set of all P-positions is identical to $S_2$, was posed by A. S. Fraenkel at the GONC 2011 workshop at the Banff Centre. Here we resolve the general case for the set of (candidate) P-positions being $S_k$. Henceforth we omit the word ``candidate'' and simply talk about sets of P-positions. We have also added the requirement that the ruleset be {\emph{invariant}} \cite{DR, LHF}, that is, the available moves do not depend on the position (for all options with non-negative coordinates). This criterion is implicitly fulfilled by many classical removal games, e.g. Nim, $k$-Wythoff Nim, Subtraction games \cite{WW} and S. Golombs take away-games \cite{Go}. Without the requirement of invariance, one may define the most trivial game rules, no move is possible from a position in $S$, and otherwise each position has a move to $(0, 0)$. On the other hand, the problem of finding invariant (but not necessarily simple) game rules for any set of P-positions, defined by a complementary pair of homogeneous Beatty sequences, was resolved in \cite{LHF}. However those game rules are not simple in the meaning that the only known formula for the invariant moves is exponentially slow in $\log(xy)$. See Figures \ref{F2}, \ref{F3} and \ref{F4} for invariant games corresponding to the CFs on page \pageref{DR}, cases $k=2$.

For many classical games, such as \emph{normal play} Nim and $k$-Wythoff Nim, the final winning position is unique, namely $(0,0)$. Given our set of P-positions, $S_k$, this requirement clearly needs to be satisfied. A convenient way to achieve this is to follow the example of our classical games, to include the Nim rules to our new game. An immediate benefit of doing this is that we automatically satisfy one of the other inherent requirements of the set $S_k$, namely that there can be at most one P-position in each row and column of $\mathbb N\times \mathbb N$. 
Precisely, the desired ruleset $\Gamma_k$ has the following permitted moves: 

\begin{theorem}\label{T1}
Let the set $S_k$ be defined by the Beatty sequences where 
$\alpha = [1;k,1,k,1,k,\ldots]$ and $\beta = [k+1;1,k,1,k,1,\ldots]$, that is 
$S_k=\{(\flrna, \flrnb)\mid n\in \mathbb N \}$. Then the invariant ruleset $\Gamma = \Gamma_k$ consisting of the following moves has a set of $P$-positions identical to the set $S_k$ (in all cases, $n,s,t \in \mathbb Z^+$):

\bigskip

\noindent {\bf{Type I - Nim Moves}}

\smallskip

$(x,y) \rightarrow (x-n,y)$ or $(x,y) \rightarrow (x,y-n)$.

\bigskip

\noindent {\bf{Type II - Extended Diagonal Moves}}

\smallskip

$(x,y) \rightarrow (x-s,y-t)$ provided that $|s-t| < k$. These moves are identical to the moves in $k$-Wythoff Nim. 

\bigskip

\noindent {\bf{Type III - Extra Moves}}

For $i = 1$ to $k-1$, use the initial value $(f^i_0,g^i_0) = (0,i+1)$ and define recursively for $n > 0$, 

$$(f^i_n,g^i_n) = (f^i_{n-1} + g^i_{n-1}, kf^i_{n-1} + (k+1)g^i_{n-1} + i)$$

The extra moves for each $i$ are $(f^i_n,g^i_n-1)$ for $n > 0$. 
\end{theorem}

Note that when $n=0$, the move $(0,i)$ is already in the ruleset as it is a Nim move. In Section \ref{S2}, we will have need to back up the recursion one step and use $(f^i_{-1},g^i_{-1}) = (-1,1)$ for all $i$.

\begin{example}\label{E4} 

For $k=2$, (and therefore $i=1$), the extra moves are
$$(2,6),(9,25),(35,96),\ldots$$
An explicit formula for the Type III moves
$(f_n^1,g_n^1)$ is given by

\begin{align}
f_n^1 &= \dfrac{(1+\sqrt 3)(2+\sqrt 3)^n + (1-\sqrt 3)(2-\sqrt 3)^n - 2}{4};\\
g_n^1 &= \dfrac{(2+\sqrt 3)^{n+1} + (2-\sqrt 3)^{n+1}-2}{2}.
\end{align}

\end{example}

Explicit formulas can be found for larger values of $k$, but are not
as succinct, and therefore we have chosen to omit them.

\begin{example} For $k=4$, the extra moves are

$$i=1:(2,10),(13,63),(77,372),\ldots$$
$$i=2: (3,16),(20,98),(119,576),\ldots$$
$$i=3:(4,22),(27,133),(163,780),\ldots$$

\end{example}

The extra moves are necessary for positions of the form $(a_n,b_n-1)$ where $c_n=1$ and $d_n = k+2$. 
For instance, when $k=4$, we seek a winning move from the N-position $(38,185)$ where $(38,186)$ is a P-position. The previous P-positions are $(37,180)$, $(36,174)$, $(35,169)$ and $(33,163)$ with differences from one P-position to its predecessor of $(1,6)$, $(1,6)$, $(1,5)$ and $(2,6)$ respectively. Preceding the nearest lesser difference of $(2,6)$ are {\emph {two}} copies of $(1,6)$ (ignoring the $(1,5)$). The winning move uses the largest valid move from the extra move set with $i=2$, namely the move $(20,98)$ which moves from $(38,185)$ to the P-position $(18,87)$.

The next section develops the machinery to examine these positions and corresponding moves. The final section shows that the rules described in Theorem \ref{T1} produce the prescribed set of P-positions.

\begin{remark}

In the case $k=1$, the set $S_1$ is the set of P-positions in 1-Wythoff Nim. The ruleset in Theorem 3 is precisely the ruleset for 1-Wythoff Nim since when $k=1$, there are no moves of Type III.

\end{remark}

\begin{remark}
In view of Figures \ref{F2} and \ref{F3}, one can see that there is in fact a very succinct description of our games as a modified greedy algorithm. Given an $S_k$-set of candidate P-positions, the algorithm starts with the moves as in $k$-Wythoff Nim as a \emph{base set of moves} and then greedily (use for example lexicographic ordering) adjoins an ordered pair of non-negative integers $(x,y)$, which does not belong to the candidate set of P-positions, as a new move if and only if the move options already defined do not suffice to find a move from $(x,y)$ to any (candidate) P-position. The new move set will be identical to our move set as in Theorem \ref{T1}. 

In this context one might want to explore other complementary Beatty sequences (forming candidate P-positions) and try and describe when similar greedy algorithms define closed formula move sets similar to the ones studied in this paper.
\end{remark}

\section{The Sturmian Word and Morphism Construction of the Beatty Sequence}\label{S2}

Here, we lay the groundwork for finding Type III winning moves for positions of the form $(a_n, b_n - 1)$, where $c_n=1$ and $d_n =k+2$. If the reader wish to become more familiar with the main structure of the proof of Theorem~\ref{T1} before reading this section, the details are given in Section \ref{S3}; all but the one most intricate case are proved without reference to Section \ref{S2}. 
Here we use some terminology from \emph{Sturmian words} and \emph{morphisms} \cite{L}. After some preliminaries, we produce the \emph{characteristic word} which corresponds to the $D$ sequence (this is Lemma \ref{lemma:3} which is proved in the Appendix) and thereby gives an alternative description of the $B$ sequence. From it, we find a new characterization of the $C$ and $A$ sequences and note some important properties. Finally, we give an algorithm for finding the desired winning move in Lemma \ref{mainlemma}. 

\begin{lemma}\label{lemma:9}
For all $n\in \mathbb{N}$ and $i\in \{1,\ldots , k-1\}$, $g^i_{n+1} = (k+2)g^i_n - g^i_{n-1}$.
\end{lemma}

\noindent{\bf {Proof:}} $g^i_{n+1} = (k+1)g^i_n + kf^i_n + i = (k+2)g^i_n - g^i_n + kf^i_n + i = (k+2)g^i_n - (kf^i_{n-1} + (k+1)g^i_{n-1} + i) + k(f^i_{n-1} + g^i_{n-1}) + i = (k+2)g^i_n - g^i_{n-1}$. \hfill $\Box$

\bigskip

\subsection{The sequence $D = (d_1,d_2,d_3,\ldots )$}

We wish to describe the sequence $D=(d_1,d_2,d_3,\ldots )$ via the Sturmian word produced by the morphism

$$\varphi(\sigma) = \sigma \tau^k = \sigma \tau \tau \tau \ldots \tau \quad (k \text{ copies of } \tau)$$

$$\varphi(\tau) = \sigma \tau^{k+1} = \sigma \tau \tau \tau \ldots \tau \quad (k+1 \text{ copies of } \tau)$$

and 

$$\varphi(uv) = \varphi(u)\varphi(v)$$ for any words $u$, $v$ consisting of the letters $\sigma$, $\tau$ where the operation is concatenation.

\begin{notation}
Let $w_0$ be the word $\sigma$, and $w_n = \varphi(w_{n-1})$. Note that $w_{n-1}$ is a prefix of $w_n$ so that $$W = \lim_{n \rightarrow \infty} \varphi^n(w_0)$$ is well-defined.
\end{notation}

\begin{example} For $k=2$, $W = \sigma \tau \tau \sigma \tau \tau \tau \sigma \tau \tau \tau  \sigma \tau \tau  \sigma \tau \tau \tau  \sigma \tau \tau \tau  \sigma \tau \tau \tau  \sigma \tau \tau \sigma \tau \tau \tau \ldots$
\end{example}

\begin{lemma}\label{lemma:3} 
For all $k$, $d_j = k+1$ if the $j^\emph{th}$ letter of $W$ is $\sigma$ and $d_j = k+2$ if the $j^\emph{th}$ letter of $W$ is $\tau$. 
\end{lemma}

We give the proof of Lemma \ref{lemma:3} in the Appendix. 

\subsection{The sequence $C = (c_1,c_2,c_3,\ldots )$}

Next, we give an equivalent construction of the sequence $C = (c_1,c_2,c_3,\ldots )$ via the $D$ sequence.

\begin{notation} The {\emph{remainder}} or {\emph{fractional part}} of $x\in \mathbb{R} $ shall be denoted $\{x\} = x - \lfloor x\rfloor $. 
\end{notation}

\begin{figure}[ht!]
  \centering
    \includegraphics[width=0.45\textwidth]{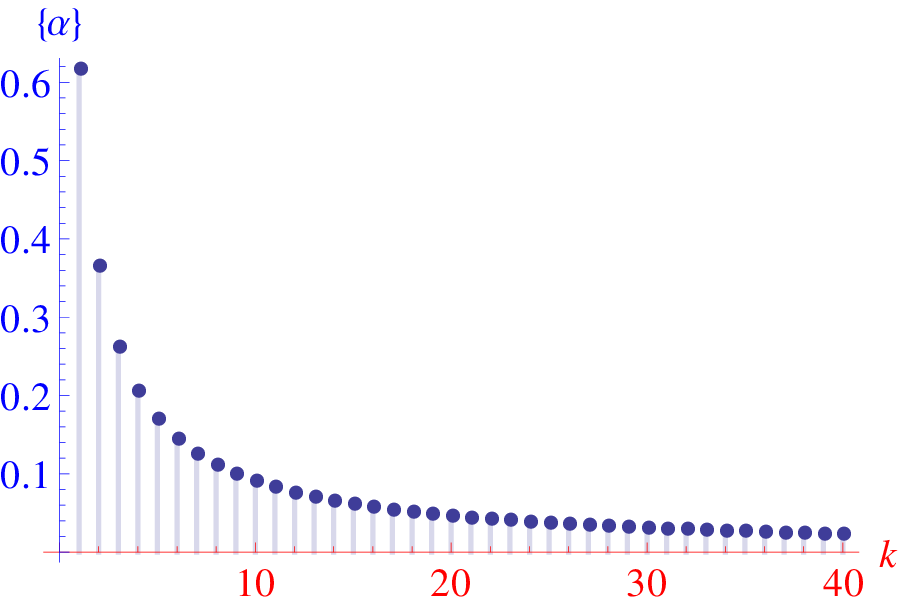} \includegraphics[width=0.45\textwidth]{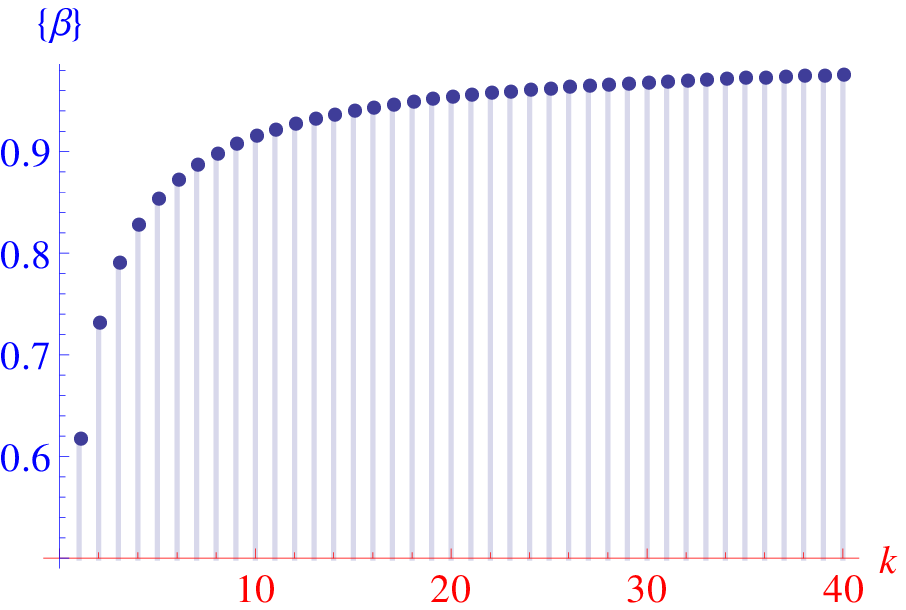} \caption{The fractional parts $\{\alpha_k\}$ and $\{\beta_k\}$, for $k\in \{1,40\}$.}\label{Ffrac}
\end{figure}

\begin{lemma}\label{lemma:10} For all $k$, $\{\alpha_k\} < 1/k$.
\end{lemma}

\noindent{\bf {Proof:}} $1 < \alpha < 2$, so $\{\alpha\} = \alpha -1 = \frac{1}{2}\left(\sqrt{1+\frac{4}{k}} - 1\right)$. $$\text{Since }(k+2)^2 > k^2\left(1 + \frac{4}{k}\right)$$

 $$\Longrightarrow \frac{k+2}{k} > \sqrt{1+\frac{4}{k}}$$ 
 
 $$\Longrightarrow \frac{1}{k} > \frac{1}{2}\left(\sqrt{1+\frac{4}{k}}-1\right) =\{\alpha\}.$$ \hfill $\Box$
 
 \begin{lemma}\label{lemma:11} For all $k$, $\{\beta_k\} = k\{\alpha_k\}$.
 \end{lemma}
 
\noindent{\bf {Proof:}} $$\beta = k\alpha + 1 = k+1 + k\{\alpha\}$$
 
 $$\Longrightarrow \{\beta\} = k\{\alpha\}$$ since $\{\alpha\} < \frac{1}{k} \Longrightarrow k\{\alpha\} < 1$. \hfill $\Box$
 
 \bigskip

We now compare the the sequences $C$ and $D$, first with an example

\begin{example} For $k=4$, 

$$\left( {C \atop D} \right)=
\left(  {1\atop 5} {1\atop 6} {1\atop 6} {1\atop 6} {2\atop 6} 
 {1\atop 5} {1\atop 6} {1\atop 6} {1\atop 6} {2\atop 6}  {1\atop 6}
  {1\atop 5} {1\atop 6} {1\atop 6} {2\atop 6}  {1\atop 6} {1\atop 6} 
   {1\atop 5}{1\atop 6} {2\atop 6} {1\atop 6} {1\atop 6} {1\atop 6}
    {1\atop 5} {2\atop 6} {1\atop 6} {1\atop 6} {1\atop 6} {2\atop 6} 
     {1\atop 5} {1\atop 6} {1\atop 6}{1\atop 6}{2\atop 6}{1\atop 5}\cdots \right)
$$

\end{example}

If we remove each 5 in the $D$ sequence and the corresponding 1 in the $C$ sequence, what remains in the $C$ sequence is periodic with the value $2$ in positions $4,8,12,\ldots$ and the value $1$ otherwise. It turns out that this observation corresponds to an alternative description of the $C$ sequence provided by the $D$ sequence, for general $k$.

\begin{lemma}\label{lemma:12} Suppose that $c_p = c_q = 2$ for $p > 0$ and $q > p$ minimal. Then there are exactly $k-1$ values of $i$, $p < i < q$ for which $d_i = k+2$.
 \end{lemma}
  
\noindent{\bf {Proof:}} Since $c_p = a_p - a_{p-1} = 2$, Lemma \ref{lemma:10} gives $\{p\alpha\} < \frac{1}{k}$. Let $i\in\{p+1,\ldots ,q-1\}$ so that $c_i = 1$ and so (the latter inequality is by Lemma \ref{lemma:10}), 
\begin{align}\label{a1k}
0<\{\alpha\} = \{i\alpha\} - \{(i-1)\alpha\}<\frac{1}{k}. 
\end{align}
By Lemma \ref{lemma:11}, this gives $\{\beta\}= k\{i\alpha\} - k\{(i-1)\alpha\}$. Now, as we have seen, going from $b_{i-1}$ to $b_i$ produces either the difference $d_i =k+1$ or $k+2$. Then, by $\beta\in (k+1,k+2)$, it is clear that the greater value will be attained if and only if there is a $j \in \{1,\ldots ,k-1\}$ such that $\{(i-1)\alpha\} < \frac{j}{k} < \{i\alpha\}$. 

By the last inequality in (\ref{a1k}), each $j$ will correspond to a unique $i$. Hence $d_i = k+2$ occurs exactly $k-1$ times between consecutive occurrences of $c_n = 2$.  \hfill $\Box$

\begin{lemma}\label{lemma:13}
Let $k\in \mathbb Z^+$. Then $c_n = 2$ iff $d_n = k+2$ and $b_n \equiv n$ (mod $k$).
\end{lemma}
\noindent{\bf Proof:} 
Suppose now that $d_i=k+1$ so that 
\begin{align}\label{b1k}
0<\{\beta\} = \{i\beta\} - \{(i-1)\beta\}<1. 
\end{align}
Then $0<k\{\alpha\} = \{ik\alpha\} - \{(i-1)k\alpha\}<1$.
If in addition $c_i = 2$ we get that $0 < \{ik\alpha\} = k\{i\alpha\} < 1$, so that $0 < \{i\alpha\} - \{\alpha\} = \{(i-1)k\alpha\}/k < 1/k$. This gives 
that $\{i\alpha\} > 1/k$ so that $c_i = 1$. Thus, we have proved that $c_n = 2$ 
implies $d_n = k + 2$. But then Lemma \ref{lemma:12} gives that $b_n-(k+1)n\equiv 0 \pmod k$, for each $n$ such that $c_n=2$. \hfill $\Box$\\

For record keeping purposes, we index the $\tau$ in the word $W$ with period $k$ so that 

$$W = \sigma \tau_1\ldots \tau_k \sigma \tau_1\ldots \tau_k \tau_1 \sigma \tau_2 \ldots$$

\begin{definition} A {\emph {syllable}} of $W$ is a string of letters of the form $\varphi(\sigma)$ or $\varphi(\tau)$, that is, it begins with $\sigma$, and ends with the $\tau$ which precedes the next $\sigma$.
\end{definition}
Thus the morphism $\varphi$ maps letters to syllables. 
Note that the indexing for the $\varphi(\sigma)$ will depend on the preceding syllable, but each index will appear exactly once. Hence, for all $i$, we get that $$\varphi(\tau_i) = \sigma \tau_i \tau_{i+1}\ldots \tau_k \tau_1\ldots \tau_i.$$ Using this notation Lemma \ref{lemma:13} states that 

\begin{align}\label{16}
c_n = 2 \text{ iff the } n^\text{th}\text{ letter of } W\text{ is }\tau_k.
\end{align}

\subsection{Sums of Factors}

\begin{definition}\label{factor} A \emph {factor} of a word is a sequence of consecutive letters. If the factor begins with the first letter of the word, the factor is called a \emph {prefix}. If the factor contains the last letter of a finite word, the factor is called a \emph {suffix}.
\end{definition}

\begin{definition}\label{initword} For each $i\in \{1,\ldots ,k-1\}$, let $w^i_0$ be the word $\sigma \tau_1\ldots \tau_i$ and $w^i_n = \varphi(w^i_{n-1})$. 
\end{definition}

\begin{lemma}\label{lemma:14} For each $i\in\{1,\ldots , k-1\}$ and all $n \ge 0$, $f^i_n$ counts the number of copies of $\tau_k$ in the word $w^i_n$ and $g^i_n$ counts the number of letters. Note that for $n \ge 1$, $g^i_{n-1}$ counts the number of syllables in the word (which equals the number of copies of $\sigma$ in $w^i_n$ by construction).
\end{lemma}

\noindent{\bf {Proof:}} {\underline {Base Case:}} $f^i_0 = 0$ and $w^i_0$ contains no $\tau_k$. $w^i_0$ contains $i+1$ letters and $g^i_0 = i+1$.

{\underline {Induction:}} The morphism $\varphi$ sends each $\tau_k$ to a syllable containing two $\tau_k$ and all other letters to a syllable containing  a single $\tau_k$, hence the number of copies of $\tau_k$ in $w^i_{n+1}$ is $2f^i_n + (g^i_n-f^i_n) = f^i_{n+1}$. The number of letters in the new word is $k+2$ for each letter subtracting one for each $\sigma$ for a total of $(k+2)g^i_n - g^i_{n-1} = g^i_{n+1}$ by Lemma \ref{lemma:9}. \hfill $\Box$

\bigskip

\begin{lemma}\label{lemma:15} A factor of $W$ of length $g^i_n$ contains either $g^i_{n-1}$ or $g^i_{n-1}-1$ copies of $\sigma$. No other number is attainable.
\end{lemma}

\noindent{\bf {Proof:}} By construction, $w^i_n$ has length $g^i_n$ and has $g^i_{n-1}$ copies of $\sigma$ so $g^i_{n-1}$ is attainable. $W$ is a Sturmian word, and therefore balanced, hence only one other value is attainable, either $g^i_{n-1} - 1$ or $g^i_{n-1} + 1$. Shift $k+2$ steps to the right in $w_n^i$. Then we lose 2 copies of $\sigma$ and gain one. Hence $g^i_{n-1} - 1$ is the correct value. \hfill $\Box$

\bigskip

\begin{lemma}\label{lemma:16} For each $i\in \{1,\ldots , k-1\}$, and for all $n \ge 1$, $(f_n^i,g_n^i)$ is a P-position.
\end{lemma}

\noindent{\bf {Proof:}} Let $j=g_{n-1}^i$, which is the length of $w_{n-1}^i$ by Lemma \ref{lemma:14}. By (\ref{16}), the number of copies of $\tau_k$ plus the number of letters in $w_{n-1}^i$ is $a_j$. By Lemma \ref{lemma:14} this equals $f^i_{n-1}+g^i_{n-1}=f^i_n$. By construction and Lemma \ref{lemma:9}, $b_j = (k+2)j - (\text{the number of copies of }\sigma \text{ in }w_{n-1}^i)=g_n^i$. \hfill $\Box$

\begin{definition}\label{def-index} A factor of $W$ has \emph{index} $i$ if it ends with $\tau_i$ for some $i~\in~\{1,\ldots , k-1\}$. A P-position $(a_n,b_n)$ has \emph{index} $i$ if the prefix of $W$ of length $n$ has index $i$.\end{definition}

\begin{lemma}\label{lemma:17}
For a fixed index $i$, let $x$ be a factor of the word $w^i_{n+2}$, with the following properties:
\begin{itemize}
\item $x$ has length $g^i_{n+1}$
\item $x$ is not the suffix of $w^i_{n+2}$
\item $x$ ends in $\tau_i$
\end{itemize}
Then $x$ contains precisely $g^i_n$ copies of $\sigma$. By construction, two equal length factors of $W$ with the same index and the same number of copies of $\sigma$ will correspond to two equal length factors of $C$ with the same number of copies of $2$. Hence the two factor sums in $C$ are equal and the two factor sums in $D$ are equal.
\end{lemma}

\noindent{\bf {Proof:}} Note that the statement is vacuous if $n=-2$

\noindent  {\underline {Base Case: $n = -1$}}

If $n = -1$, then $x$ has length $i+1$. $w_1^i$ has $i+1$ syllables, with $\tau_i$ in position $i+1$ in the first syllable and in position $i+3-s$ in syllable $s$ for $2 \le s \le i+1$, hence the $i+1$ letters ending in $\tau_i$ always contain exactly one $\sigma$.

\noindent  {\underline {Induction:}}

If the terminal $\tau_i$ of the factor $x$ is the last letter of a non-terminal syllable of $w^i_{n+2}$, then the factor contains exactly $g^i_n$ syllables since the terminal $\tau_i$ was a result of the output of $\varphi(\tau_i)$, and the previous word $w^i_{n+1}$ has precisely $g^i_{n-1}$ copies of $\sigma$ in a factor of length $g^i_n$ by induction. 

If the terminal $\tau_i$ is not the last letter in its syllable, then compare the factor $x$ with the nearest previous factor $y$ for which the terminal $\tau_i$ is the last letter in its syllable. If the factors $x$ and $y$ overlap so that there exist non-empty words $t$, $u$, $v$ with $y = tu$, $x=uv$, we need to show that the number of copies of $\sigma$ in $t$ equals the number of copies of $\sigma$ in $v$. 

Let $j$ be the index of the syllable containing the terminal $\tau_i$ of the factor $x$. If there is no syllable $\varphi(\sigma)$ in $v$, then the length of $v$ is $(k+1)m$ where $m = j-i$ if $j>i$ and $m=k+j-i$ if $j \le i$. $v$ contains $m-1$ full syllables plus the terminal partial syllable. Each full and partial syllable contains one $\sigma$, so $v$ contains $m$ copies of $\sigma$. In other words, the fraction of letters in $v$ which are $\sigma$ is $\dfrac{1}{k+1}$. If $v$ does contain a syllable $\varphi(\sigma)$, this ratio is unchanged. 

For any integer $m$, the number of copies of $\sigma$ in any factor of length $(k+1)m$ cannot exceed this ratio since the length of each syllable is $\ge k+1$. Since the number of copies of $\sigma$ in $y$ is $g_n^i
$ which is maximal by Lemma \ref{lemma:15}, the number of copies of $\sigma$ in $x$ which is at least as many as in $y$ must also be $g_n^i$ hence the number of copies of $\sigma$ in $y$ equals the number of copies of $\sigma$ in $x$.

In the case that $x$ and $y$ do not overlap, note that the maximum
distance that $x$ needs to be shifted occurs when the terminal
$\tau_i$ of $x$ is the leading $\tau$ in a syllable ending in $\tau_i$
and that this distance is $k-1$ syllables of length $k+2$ plus perhaps
a syllable of length $k+1$ plus 2 for a total of $[(k-1)(k+2)]+[k+1]+2
= (k+1)^2 < g_2^i$ thus in the induction step, $x$ and $y$ do not
overlap only for $n=0$. In this case, $y$ has one syllable of length
$k+1$ and $i$ syllables of length $k+2$ yielding $i+1$ copies of
$\sigma$. $x$ has its terminal partial syllable of length $\le k+1$,
thus $x$ has at least as many copies of $\sigma$ as does $y$, and
since the number of copies of $\sigma$ in $y$ is maximal by Lemma
\ref{lemma:15}, the number of copies of $\sigma$ is the same in $x$
and $y$. \hfill $\Box$\\

At the beginning of this section we promised an algorithm for finding a certain winning move. We deliver it here:

\begin{lemma}\label{mainlemma} Let $(x,y) = (a_n,b_n)$ be a P-position
with index $i \in \{1,\ldots,k-1\}$. From the position $(x,y-1)$, the
Type III move $(u,v)$ corresponding to $i$ with $v \le b_n - 1$
maximal is to a P-position.
\end{lemma}

{\bf {Proof:}} Find $m$ such that $g^i_m \le b_n < g^i_{m+1}$. In the
first case, if $b_n = g^i_m$, then from $(x,y-1)$, the extra move
$(f^i_m,g^i_m - 1)$ is to $(0,0)$ by Lemma \ref{lemma:16}. In all
other cases, Lemma \ref{lemma:17} shows that all factors with index
$i$ and length $g_{m-1}^i$, except the last, in the word $w^i_{m}$
have the same number of copies of $\sigma$, hence the factor sums in
$C$ and $D$ have the same sums as in the first case, namely $f^i_m$
and $g^i_m - 1$. Hence the move $(f^i_m,g^i_m-1)$ is to the P-position
$(a_{n-j},b_{n-j})$ where $j = g^i_{m-1}$.  \hfill $\Box$

\section{The rules are correctly defined}\label{S3}
In this section, we prove Theorem \ref{T1}, that is we verify that the set $S_k$ is generated as the complete set of P-positions by the ruleset $\Gamma_k$.
\begin{definition} The {\emph{gap}} of a P-position $(a_n,b_n)$, denoted $\delta_n$ is
$$\delta_n = b_n-a_n$$
The {\emph{gap difference}} between two P-positions $(a_m,b_m)$ and $(a_n,b_n)$ with $m>n$, denoted $\Delta(m,n)$ is 
$$\Delta(m,n) = \delta_m - \delta_n$$
\end{definition}

\bigskip

We must check that there is no move connecting any two  
P-positions (such a ``short-circuit'' would force at least one 
of the P-positions to be de facto N and so we had to exclude 
it from the set $S_k$) and that every N-position has a move 
to a position in the set $S_k$ (for otherwise one of the 
N-positions would be de facto P, and so we had to include it to 
the set $S_k$).

\bigskip

\noindent{\bf {Proof - Part I - No move connects P to P:}} 

\bigskip

By the complementarity of the Beatty sequences, moves of Type I cannot connect any two P-positions. 

\bigskip

Note that $\Delta(m,m-1) = k$ or $\Delta(m,m-1) = k+1$. Recall that $a_m - a_{m-1} \le 2$, so moves of Type II cannot 
connect $(a_m,b_m)$ and $(a_{m-1},b_{m-1})$. If $m-n>1$, then $\Delta(m,n) \ge 2k$, so moves of Type II cannot connect P-positions in this case either.

\bigskip

It remains to justify that moves of Type III never connect two P-positions. Let $(p,q)$, $p<q$ be an extra move so that $(p,q+1) = (a_i,b_i)$, for some positive integer $i$, by Lemma \ref{lemma:16}. From the P-position $(a_m,b_m)$, it is clear that $(a_m-q,b_m-p)$ is not a P-position since $(b_m-p)-(a_m-q) > \delta_m$ and the gap must decrease.

To show that $(a_n-p,b_n-q)$ is not a P-position, assume the contrary, and note that 
\begin{align}\label{something0}
(a_n-p,b_n-q)  &= (a_n-a_i, b_n - b_i +1)\\ 
&= (\flrna-\flria,\flrnb-\flrib+1).\notag
\end{align} 
We have
\begin{align*}
\flrnia &= (n-i)\alpha - \{(n-i)\alpha \}\\&= n\alpha - i\alpha -\{(n-i)\alpha \}\\ &= \flrna + \{n\alpha\} - \flria - \{i\alpha\} - \{(n-i)\alpha\},
\end{align*}
but then 
\begin{align}\label{something}
\lfloor (n-i)\alpha \rfloor - \flrna + \flria = \{n\alpha \} - \{i\alpha\} - \{(n-i)\alpha\}
\end{align}
must be an integer and is therefore either $0$ (if  $\{n\alpha \} \ge \{i\alpha\} + \{(n-i)\alpha\}$) or $-1$ (if  $\{n\alpha \} < \{i\alpha\} + \{(n-i)\alpha\}$).\\

\noindent {\sc Case 1:} $\{n\alpha\} - \{i\alpha\} - \{(n-i)\alpha\} = 0$. Then, (\ref{something0}) and (\ref{something}) give that  
\begin{align*}
a_n-p &= \flrna - \flria\\ &= \flrnia\\ &= a_{n-i}. 
\end{align*}
Hence, for $(a_n-p,b_n-q)$ to be a $P$-position, we must have 
\begin{align*}
b_n - b_i + 1 &= b_{n-i}\\ 
&= \flrnib\\
&=\lfloor b_n + \{n\beta\} - b_i - \{i\beta\} \rfloor\\ 
&= b_n - b_i + \lfloor \{n\beta\}-\{i\beta\} \rfloor,
\end{align*} 
but $\lfloor \{n\beta\}-\{i\beta\} \rfloor$ cannot be $1$.\\

\noindent{\sc Case 2:} $\{n\alpha\} - \{i\alpha\} - \{(n-i)\alpha\} = -1$. Then (\ref{something}) gives that 
$$\flrna - \flria = \flrnia + 1.$$ By the latter expression, this number, which is strictly greater than zero, can belong either to the set $A$ or $B$.
If  $$\flrnia + 1 = b_x\in B,$$ then, for $$(a_n-p,b_n-q)=(\flrna-\flria,\flrnb-\flrib +1)$$ to be a non-trivial P-position, we must have that $a_x = b_n-q$ and $b_x = a_n-p$. But, since for all $x>0$, $b_x>a_x$, this gives
\begin{align*}
a_n-a_i &> b_n - b_i +1, 
\end{align*}
which is false, since $\delta_n>\delta_i$ if $n>i$. 

Otherwise 
\begin{align*}
\flrna - \flria &=a_{n-i+1}\in A, 
\end{align*}
so that, by (\ref{something0}) and the definition of a P-position, we must have 
\begin{align}\label{something1}\notag
b_n - b_i + 1 &= b_{n-i+1}\\ 
&= \lfloor (n-i+1)\beta \rfloor. 
\end{align}
However 
\begin{align}\label{something2}\notag
\lfloor (n-i+1)\beta \rfloor &= \lfloor b_n + \{n\beta\} - b_i - \{i\beta\} +b_1 - \{\beta\} \rfloor\\ &= b_n - b_i + b_1 + \lfloor \{n\beta\} - \{b_i\} - \{b_1\} \rfloor. 
\end{align}
The last term is either $0$ or $-1$. There are no moves of Type III for the case $k=1$, thus $k \ge 2$ and  $\beta > 3$. Therefore $b_1 \ge 3$, which gives $b_1 + \lfloor \{n\beta\} - \{i\beta\} - \{\beta\} \rfloor \ne 1$, which, by (\ref{something2}), contradicts (\ref{something1}).

\bigskip

\noindent{\bf {Proof - Part II - Every N has a move to a P:}} 

Assume in all cases that $x \le y$.

If $(x,y)$ is an N-position and either $x \in B$ or $y \in B$, then there is a Nim move (Type I) to a P-position. If $x = a_n \in A, y \in A$, with $y > b_n$, the Nim move lowering $y$ to $b_n$ is to a P-position.

If $x = a_n, y \in A, y < b_n -1$, then $y-x \le \delta_n -2$. Since the gaps $\delta_j$ increase by either $k$ or $k+1$ as $j$ increases by 1, then there is an Extended Diagonal move (Type II) to a P-position corresponding to $\delta_j$ which is nearest $y-x$.

What remains to be shown are winning moves from $x =a_n, y = b_n -1$. If $a_n = a_{n-1} + 2$ or $b_n = b_{n-1} + k+1$, then the extended diagonal move $(2,k+1)$ or $(1,k)$ moves to the P-position $(a_{n-1},b_{n-1})$. Otherwise, Lemma \ref{mainlemma} finds the winning Type III move.
 \hfill $\Box$\\

\appendix
\numberwithin{equation}{section}
\numberwithin{theorem}{section}
\setcounter{theorem}{0}

All notation in the appendix is local unless stated otherwise. We use theory from \cite{Gl} (further references are given in \cite{Gl}). The words are defined on the alphabet 
$\{0,1\}.$ For $k\ge 2$ an integer, we are interested in the morphism
\begin{align}
\theta:\ \ & 0\rightarrow 1^k0 \label{theta0}\\ 
& 1\rightarrow 1^k01 \label{theta1}
\end{align}
which we will show corresponds to the positive root, 
\begin{align}\label{gamma}
\gamma = \frac{\sqrt{k^2+4k}-k}{2}\in (1/2,1),
\end{align}
 of $x^2+kx-k = 0.$ Namely, by a result in \cite{Gl}, 
we will obtain that $$\lim_{n\rightarrow \infty} \theta^n(1)$$ is the characteristic word $c_\gamma$ of $\gamma$. The density of (the 1s in) $c_\gamma$ is of course $\gamma$. Also $\beta $ as defined in Section \ref{S1} equals $\gamma +k+1$ (that is $\gamma = \{\beta\}$). Hence the continued fraction expansion of $\gamma$ is $\gamma =[0;1,k,\overline{1,k}]$ (where $\overline x$ denotes the periodic pattern $x,x,\ldots$).

Let $X = \lim_{n\rightarrow\infty}\theta^n(1)$.
We show by induction that the Sturmian word $W=\lim_{\rightarrow \infty}\varphi (\sigma)$, defined as in Section \ref{S2}, is identical (via $\sigma \leftrightarrow 0, \tau\leftrightarrow 1$) to the word $0X$. That is, we want to show that: 
\begin{lemma}
The $i^\emph{th} \text{ letter of } 0X \text{ is a } 1$ if and only if the $i^\emph{th} \text{ letter of } W \text{ is a } \tau$. 
\end{lemma}
\noindent {\bf Proof:}
The first letter in $0X$ and $W$ is $0$ and $\sigma$ respectively; $\varphi$ acts on its letter, whereas $\theta$ does not. Rather, $\theta$ acts on the first letter in $X$. 
 Let us state our induction hypothesis:\\ 

\noindent Case 1, $x_j=0$: Then the last letter of the $j^\text{th}$ syllable of $\theta$, as in the right hand side of (\ref{theta0}), corresponds precisely to the first letter of the $(j+1)^\text{st}$ syllable of $\varphi$.\\

\noindent Case 2, $x_j=1$: Then the last two letters of the $j^\text{th}$ syllable of $\theta$, as in the right hand side of (\ref{theta1}), correspond precisely to the first two letters of the $(j+1)^\text{st}$ syllable of $\varphi$.\\

If these two cases hold for all $j$, then, by \begin{align*}
(\text{the lenght of }\varphi(\sigma))=(\text{the length of }\theta(0))=k+1 \intertext{and} (\text{the length of }\varphi(\tau))=(\text{the length of }\theta(1))=k+2,
\end{align*}
 the infinite words correspond precisely. This follows since then the first $k$ letters in the $j^\text{th}$ syllable of $X$, each a copy of 1, correspond precisely to the last $k$ letters of the $j^\text{th}$ syllable of $W$, each a copy of $\tau$.

Our base case is that the first syllable of $X$ ends with $01$ and the second syllable of $W$ begins with $\sigma\tau$, and indeed it holds for the  prefixes $01^k01$ and $\sigma\tau^k\sigma\tau^{k+1}$ respectively. 

But then, comparing the definitions of $\varphi$ and $\theta$ with the paragraph after Case 2, the induction hypothesis gives the claim. \hfill $\Box$

Define on $\{0,1\}$ the following three morphisms
\[
  E: \begin{matrix}
      &0 &\mapsto &1& \\
      &1 &\mapsto &0&
     \end{matrix}, \qquad \eta: \begin{matrix}
                                   &0 &\mapsto &01& \\
                                   &1 &\mapsto &0& 
                                   \end{matrix}, \qquad 
  \overline{\eta}: \begin{matrix}
                  &0 &\mapsto &10& \\
                  &1 &\mapsto &0~&
                 \end{matrix}.
\]

A morphism $\psi$ is \emph{Sturmian} if and only if it is a
composition of $E$, $\eta$, and $\overline{\eta}$ in any number
and order. Furthermore, 
a morphism $\psi$ is \emph{standard} if and only if it is a
composition of $E$ and $\eta$ in some order. A morphism is 
\emph{non-trivial} if it is neither $E$ nor the identity morphism.

Suppose $\alpha = [0;1+d_1,d_2,d_3, \ldots]$, with $d_1 \geq 0$
and all other $d_n > 0$. To the \emph{directive sequence}
$(d_1,d_2,d_3,\ldots)$, we associate a sequence $(s_n)_{n \geq
-1}$ of words defined by
\[
  s_{-1} = 1, ~s_{0} = 0, ~s_{n} = s_{n-1}^{d_{n}}s_{n-2}; \quad n \geq 1.
\]
Such a sequence of words is called a \emph{standard sequence}.

For any $n\geq0$, $s_n$ is a
prefix of $s_{n+1}$, so that $\underset{n\rightarrow \infty}{\mbox{lim}}s_n$ is well defined as an infinite word.
Moreover, standard sequences are related to characteristic Sturmian words. Each $s_n$ is a prefix of $c_\alpha$, and we have
\[ 
c_{\alpha} = \underset{n \rightarrow \infty}{\mbox{lim}}s_n.
\]

In \cite{Gl}, all irrationals $\alpha \in (0,1)$ such that the
characteristic Sturmian word $c_\alpha$ is generated by a
morphism are classified. A \emph{Sturm number} is an irrational number $\alpha \in (0,1)$ that has a continued fraction expansion of one of the following types:
\begin{itemize}
\item[(i)] ~$\alpha = [0;1 + d_1,\overline{d_2,\ldots,d_n}] <
\frac{1}{2}$ with ~$d_n \geq d_1 \geq 1$;
\item[(ii)] ~$\alpha = [0;1,d_1,\overline{d_2,\ldots,d_n}] >
\frac{1}{2}$ with ~$d_n \geq d_1$.
\end{itemize}
Observe that if $\alpha = [0;1+d_1,\overline{d_2,\ldots,d_n}]$
with $d_n \geq d_1 \geq 1$, then
\[
  1 - \alpha = [0;1,d_1,\overline{d_2,
                        \ldots,d_n}].
\]
Hence, $\alpha$ has an expansion of type (i) if and only if $1 -
\alpha$ has an expansion of type (ii). Accordingly, $\alpha$ is
a Sturm number if and only if $1-\alpha$ is a Sturm number and one can show that 
$c_{1-\alpha}$ is obtained from $c_{\alpha}$ by exchanging all
letters $0$ and $1$ in $c_\alpha$, so that 

\begin{align}\label{Ec}
c_{1-\alpha} =E(c_{\alpha}).
\end{align}

Therefore, we can restrict our attention to
characteristic Sturmian words $c_\alpha$ such that $\alpha$ is a 
Sturm number of type (i).

We say that a morphism $\psi$ \emph{fixes} an infinite word $x$ if
$\psi(x) = x$, in which case $x$ is called a \emph{fixed point} of
$\psi$. The following result describes all irrationals $\alpha \in (0,1)$ such
that $c_\alpha$ is a fixed point of a non-trivial morphism. 

\begin{theorem}\cite{BS}\label{13}
 Let $\alpha \in (0,1)$ be irrational. Then 
$c_\alpha$ is a fixed point of a non-trivial morphism $\sigma$ if
and only if $\alpha$ is a Sturm number. In particular, if $\alpha
= [0;1 + d_1,\overline{d_2, \ldots, d_n}]$ with ~$d_n
\geq d_1 \geq 1$, then $c_\alpha$ is the fixed point of any power
of the morphism
\[
  \sigma: \begin{matrix}
           &0& &\mapsto &s_{n-1}~~~~~~~& \\
           &1& &\mapsto &s_{n-1}^{d_n - d_1}s_{n-2}&
          \end{matrix}.
\]
\end{theorem}

Our irrational $\gamma = [0;1,k,\overline{1,k}]$ (as in (\ref{gamma})) is of type (ii) for all $k$, with $n = 3$, $d_1 = k, d_2 = 1, d_3 = k$, and so we rather apply the theorem to $\alpha = 1-\gamma =  [0;1+k,\overline{1,k}]$ which is of type (i). For our application, we have that $s_{-1}=1, s_0 = 0, s_1=0^k1$ and $s_2=0^k10$, so that the morphism $\sigma$ in Theorem \ref{13} corresponds to $0\rightarrow 0^k10$ and $1\rightarrow 0^k1$. By (\ref{Ec}), it is easy to check that the standard morphism $(E\eta)^k\eta$ is identical to $E\sigma = \theta$, so that $E(\lim_{n\rightarrow\infty}\sigma^n(0))$ corresponds to the characteristic word $c_\gamma$ with $\gamma$ as in (\ref{gamma}). This proves Lemma \ref{lemma:3}. \hfill$\Box$

\enddocument